\documentclass[11pt]{article}
\usepackage{mathrsfs}
\usepackage{amsmath}
\usepackage{amssymb}
\usepackage{graphicx}
\def \la{\lambda}
\def \lamin{\la_{\min}}

\def \U{\mathbf{U}}

\def \I{\mathbf{I}}
\def \iff{\text{if and only if }}
\newtheorem{theorem}{\scshape \mdseries  Theorem}[section]
\newtheorem{lemma}[theorem]{\scshape \mdseries  Lemma}

\begin{document}

\title{\bf The least eigenvalue of graphs whose complements are unicyclic\thanks{
Supported by National Natural Science Foundation of China (11071002),
Program for New Century Excellent Talents in University (NCET-10-0001),
Key Project of Chinese Ministry of Education (210091),
Specialized Research Fund for the Doctoral Program of Higher Education (20103401110002),
Science and Technological Fund of Anhui Province for Outstanding Youth (10040606Y33),
Scientific Research Fund for Fostering Distinguished Young Scholars of Anhui University(KJJQ1001),
Academic Innovation Team of Anhui University Project (KJTD001B). }}
\author{Yi Wang$^1$, Yi-Zheng Fan$^{1,}$\thanks{Corresponding author. E-mail addresses: wangy@ahu.edu.cn (Y. Wang),
fanyz@ahu.edu.cn (Y.-Z. Fan), starlittlestar@sina.com (X.-X. Li),
zhangfeifei2403@126.com (F.-F. Zhang). }
, Xiao-Xin Li$^2$, Fei-Fei Zhang$^1$
\\
  {\small  \it $^1$School of Mathematical Sciences, Anhui University, Hefei 230039, P. R. China} \\
  {\small \it $^2$Department of Mathematics and Computer Sciences, Chizhou University, Chizhou 247000, P.R. China}
 }
\date{}
\maketitle

{\small
\noindent{\bf Abstract:} A graph in a certain graph class is called minimizing 
if the least eigenvalue of the adjacency matrix of the graph attains the minimum among all graphs in that class.
Bell {\it et al.} have characterized the minimizing graphs in the class of connected graphs of order $n$ and size $m$, whose 
complements are either disconnected or contain a clique of order at least $n/2$.
In this paper we discuss the minimizing graphs of a special class of graphs of order $n$ whose complements are connected and contains exactly one cycle 
(namely the the class $\mathscr {U}^{c}_{n}$ of graphs whose complements are unicyclic), 
and characterize the unique minimizing graph in $\mathscr {U}^{c}_{n}$ when $n \geq 20$.


}

\section{\bf Introduction}
Let $G=(V,E)$ be a simple graph with vertex set $V=V(G)=\{v_1,v_2,\ldots,v_n\}$ and edge set $E=E(G)$.
The {\it adjacency matrix} of $G$ is defined to be a matrix $A(G)=[a_{ij}]$ of order $n$,
  where $a_{ij}=1$ if $v_i$ is adjacent to $v_j$, and $a_{ij}=0$ otherwise.
Since $A(G)$ is real and symmetric, its eigenvalues are real and can be arranged as:
   $\la_1(G) \leq \la_2(G) \leq \cdots \leq \la_n(G)$.
We simply call the eigenvalues of $A(G)$ as {\it the eigenvalues} of $G$.
The eigenvalue $\la_n(G)$ is exactly the spectral radius of $A(G)$;
and there are many results in literatures concerning this eigenvalue; see e.g. \cite{bru,cve1,cve2}.

The least eigenvalue $\la_1(G)$ is now denoted by $\lamin (G)$,
   and the corresponding eigenvectors are called the {\it first eigenvectors} of $G$.
Relative to the spectral radius, the least eigenvalue has received less attention.
In the past the main work on the least eigenvalue of a graph is focused on its bounds; see e.g. \cite{cve3,hong}.
Recently, the problem of  minimizing the least eigenvalues of graphs subject to graph parameters
   has received much more attention, since two papers of Bell {\it et al.} \cite{bell1, bell2} and one paper
  of our group \cite{fan} appeared in the same issue of the journal Linear Algebra and Its Applications.
Ye and Fan \cite{ye} discuss the connectivity and the least eigenvalue of a graph.
Liu {\it et al.} \cite{liu} discuss the least eigenvalues of unicyclic graphs with given number of pendant vertices.
Petrovi\'{c} {\it et al.} \cite{pe} discuss the least eigenvalues of bicyclic graphs.
Wang {\it et al.} \cite{wang,wang2} discuss the least eigenvalue and the number of cut vertices of a graph.
Tan and Fan \cite{tanfan} discuss the least eigenvalue and the vertex/edge independence number, the vertex/edge cover number of a graph.

For convenience,  a graph is called {\it minimizing } in a certain graph class
  if its least eigenvalue attains the minimum among all graphs in the class.
Let $\mathscr{G}(n,m)$ denote the class of connected graphs of order $n$ and size $m$.
Bell {\it et al.} \cite{bell1} have characterized the structure of the minimizing graphs in $\mathscr{G}(n,m)$; see Theorem \ref{bell}.

\begin{theorem}{\em \cite{bell1}}\label{bell}
Let $G$ be a minimizing graph in $\mathscr{G}(n,m)$.
Then $G$ is either

{\em (i)} a bipartite graph, or

{\em (ii)} a join of two nested split graphs (not both totally disconnected).
\end{theorem}

\noindent
We find that the complements of the minimizing graphs in $\mathscr{G}(n,m)$ are either disconnected or contain a clique of order at least $n/2$.
This motivates us to discuss the least eigenvalue of graphs whose complements are connected and contain clique of small size.
In a recent work \cite{fan1} we characterize the unique minimizing graph in the class of graphs of order $n$ whose complements are trees.

In this paper, we continue this work on the complements of unicyclic graphs,
  and determine the unique minimizing graph in $\mathscr {U}^{c}_{n}$ for $n \geq 20$,
  where $\mathscr {U}^{c}_{n}$ denotes the class of the complements of connected unicyclic graphs of order $n$.
It is easily seen that $\mathscr {U}^{c}_{n} \subsetneq \mathscr{G}\left(n,{n \choose 2}-n\right)$.
However the minimizing graph in $\mathscr {U}^{c}_{n}$ does not hold the conditions in Theorem \ref{bell}.

\section{\bf Preliminaries}
We begin with some definitions.
Given a graph $G$ of order $n$, a vector $X \in \mathbb{R}^n$ is called to be defined on $G$,
  if there is a 1-1 map $\varphi$ from $V(G)$ to the entries of $X$; simply written
  $X_u=\varphi(u)$ for each $u \in V(G)$.
If $X$ is an eigenvector of $A(G)$, then it is naturally defined on $V(G)$, {\it i.e.}
$X_u$ is the entry of $X$ corresponding to the vertex $u$.
One can find that
   $$X^{T}AX = 2 \sum_{uv \in E(G)} X_u X_v,\eqno(2.1)$$
  and $\la$ is an eigenvalue of $G$ corresponding to the eigenvector $X$ if and only if $X \neq 0$ and
   $$\la X_v = \sum_{u \in N_G(v)} X_u,  \hbox{~ for each vertex ~} v \in V(G).  \eqno(2.2)$$
where $N_G(v)$ denotes the neighborhood of $v$ in the graph $G$.
The equation (2.2) is called {\it $(\la,X)$-eigenequation} of $G$.
In addition, for an arbitrary unit vector $X \in \mathbb{R}^{n}$,
   $$\lamin(G) \le X^{T}A(G)X, \eqno(2.3)$$
with equality if and only if $X$ is a first eigenvector of $G$.

In this paper all unicyclic graphs are assumed to be connected.
Denote by $\mathscr{U}_n$ the set of unicyclic graphs of order $n$, and denote
$\mathscr{U}^c_n=\{G^c: G \in \mathscr{U}_n\}$, where $G^{c}$ denotes the complement of the graph $G$.
Note that $A(G^{c}) = \mathbf{J} - \mathbf{I} - A(G)$, where $\mathbf{J},\mathbf{I}$ denote
 the all-ones matrix and the identity matrix both of suitable sizes, respectively.
So for any  vector $X$,
    $$X^{T}A(G^{c})X = X^{T}(\mathbf{J} - \mathbf{I} )X - X^{T}A(G)X.   \eqno(2.4)$$

A star of order $n$, denoted by $K_{1,n-1}$, is a tree of order $n$ with $n-1$ pendant edges attached to a fixed vertex.
The vertex of degree $n-1$ in $K_{1,n-1}$ is called the {\it center} of $K_{1,n}$.
A cycle and a complete graph both of order $n$ are denoted by $C_n,K_n$ respectively.
Denote by $S_n^3$  the graph obtained from $K_{1,n-1}$ by adding a new edge between two pendant vertices.
Next, we introduce two special unicyclic graphs denoted by $\U(p,q)$ and $\U'(p)$ respectively,
where $\U(p,q)$ is obtained from two disjoint graphs $K_{1,p}\;(p \ge 1)$ and $S_{q+1}^3 (q \ge 3)$
   by adding a new edge between one pendant vertex of $K_{1,p}$ and one pendant vertex of $S_{q+1}^3$,
   $\U'(p)$ is  obtained from two disjoint graphs $K_{1,p}\;(p \ge 1)$ and $C_3$
   by adding a new edge between one pendant vertex of $K_{1,p}$ and one  vertex of $C_3$; see Fig. 2.1.

\begin{center}
\vspace{3mm}
\includegraphics[scale=.68]{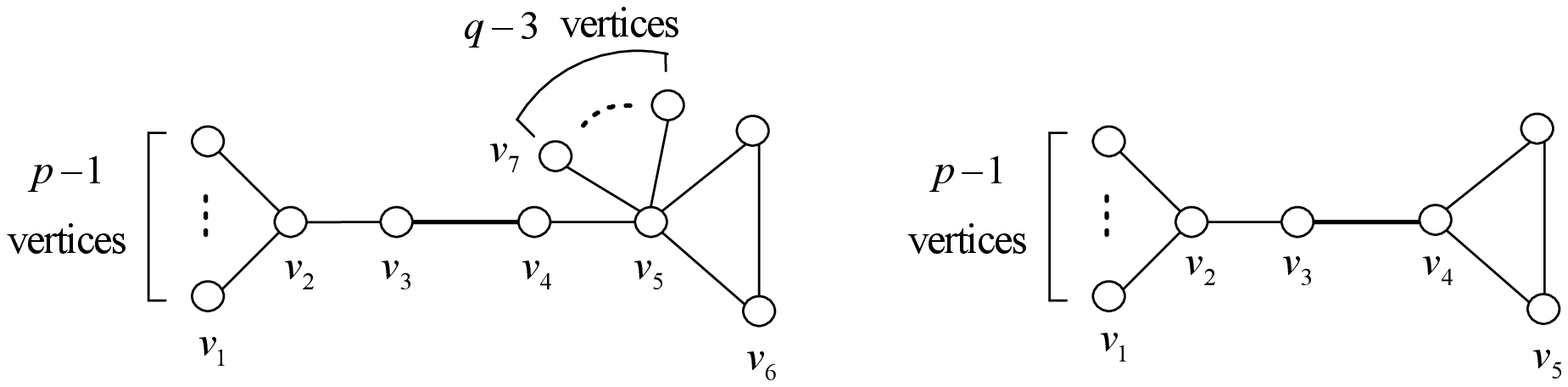}
\vspace{2mm}

{\small Fig. 2.1. The graphs $\U(p,q)$ (left side) and $\U'(p)$ (right side)}

\end{center}

For a graph $G$ containing at least one edge, then $\lamin(G) \le -1$ with equality if and only if
  $G$ is a complete graph or a union of disjoint complete graphs at least one component of which is nontrivial (or contains more than one vertices).
So, for a unicyclic graph $U$ other than $C_4$, $\lamin(U^c)<-1$ .
In addition, if $U^c$ is disconnected, then $U$ contains a complete multipartite graph as a spanning subgraph,
  which implies $U$ is $C_4$ or $S_n^3$.
When $n \ge 4$,  $(S_n^3)^c$ consists of an isolated vertex and a connected non-complete subgraph of order $n-1$.

At the end of this section, we will discuss the least eigenvalues of $\U(p,q)^c$ and $\U'(p)^c$.
Let $X$ be a first eigenvector of the graph $\U(p,q)^{c}$ with some vertices labeled in Fig. 2.1.
By eigenequations (2.2), as $\lamin(\U(p,q)^c)<-1$,
   all the pendant vertices attached at $v_2$ have the same value as $v_1$ given by $X$, say $X_1$.
Similarly all the pendant vertices attached at $v_5$ have the same value as $v_7$, say $X_7$;
  two vertices of degree $2$ on the triangle have the same value as $v_6$, say $X_6$.
Write $X_{v_i}=:X_i$  for $i=2,3,4,5$ and $\lamin(\U(p,q)^c)=:\la_1$ simply.
Then by eigenequations (2.2) on $v_i$ for $i=1,2,\ldots,7$, we have
   $$\left\{
  \begin{array}{lcl}
   \la_1 X_{1} & = & (p - 2)X_{1} + X_{3} + X_{4} + X_{5} + 2X_{6}+(q-3) X_7 \\
   \la_1 X_{2} & = & X_{4} + X_{5} +2X_{6}+ (q - 3)X_7\\
   \la_1 X_{3} & = & (p - 1)X_{1} + X_{5} +2X_{6}+(q - 3)X_7\\
   \la_1 X_{4} & = & (p - 1)X_{1} + X_{2} +2X_{6}+ (q - 3)X_7 \\
   \la_1 X_{5} & = & (p - 1)X_{1} + X_{2} + X_{3}\\
   \la_1 X_{6} &= & (p- 1)X_{1} + X_{2} + X_{3} + X_{4} + (q - 3)X_7\\
   \la_1 X_{7} &= & (p- 1)X_{1} + X_{2} + X_{3} + X_{4} + 2X_{6}+(q - 3)X_7
   \end{array}
   \right. \eqno(2.5)$$
Transform (2.5) into a matrix equality $(B-\la_1 \I)X'=0$,
  where $X'=(X_1,X_2,\ldots,X_7)^T$ and the matrix $B$ of order $7$ is easily to seen.
We have
{\small
   $$\begin{array}{ccl}
   f(\la;p,q)& : = &  \det (B-\la \I)\\
   & = &  (-8+2p+2q) + (13-11p-7q+4pq)\la +
  (20-6q-4qp)\la^2 + (-1+11p+7q-7pq)\la^3\\
   & & +  (-20+12p+12q-2pq)\la^4 + (-16+6p+6q)\la^5 + (-6+p+q)\la^6 - \la^7.
   \end{array} \eqno(2.6)$$
   }
It is easily found that $\la_1$ is the least root of  the polynomial $f(\la;p,q)$.

Let $Y$ be a first eigenvector of the graph $\U'(p)^{c}$ with some vertices labeled in Fig. 2.1.
By a similar discussion, all the pendant vertices attached at $u_2$ have the same values given by $Y$, say $Y_1$.
Two vertices of degree 2 on the triangle have the same values, say $Y_5$.
Write $Y_{u_i}=:Y_i$  for $i=2,3,4$ and $\lamin(\U'(p)^c)=:\la'_1$ simply.
Then by eigenequations on $u_i$ for $i=1,2,\ldots,5$,
  $$\left\{
  \begin{array}{lcl}
   \la'_1 Y_{1} & = & (p - 2)Y_{1} + Y_{3} + Y_{4} + 2Y_{5} \\
  \la'_1 Y_{2} & = & Y_{4} + 2Y_{5} \\
  \la'_1 Y_{3} & = & (p - 1)Y_{1} + Y_{4}+ 2Y_{5}\\
  \la'_1 Y_{4} & = & (p - 1)Y_{1} + Y_{2}\\
  \la'_1 Y_{5} & = & (p - 1)Y_{1} + Y_{2} + Y_{3}\\
  \end{array}
  \right. \eqno(2.7) $$
It is easily found that $\la'_1$ is the least root of the following polynomial:
  $$ g(\la;p): = (-4+2p) + (3-5p)\la +(6-p)\la^2 + (1+4p)\la^3 +  (-2+p)\la^4 -\la^5 .\eqno(2.8) $$

\begin{lemma}
If $n \geq 13$, $\lamin(\U(n-5,3)^c) < \lamin(\U'(n-4)^c)$.
\end{lemma}

{\bf Proof:}
Write $\lamin(\U(n-5,3)^c)=:\la_1$, $\lamin(\U'(n-4)^c)=:\la'_1$ simply.
By above discussion, $\la_1$ (respectively, $\la'_1$) is the least root of $f(\la;n-5,3)$ (respectively, $g(\la;n-4)$).
Denote by $$\bar{g}(\la;n-4):=(\la+1)^2g(\la;n-4).$$
Since $\la'_1 <-1$, $\la'_1 $ is also the least root of  $\bar{g}(\la;n-4)$.
From (2.8),  $g(-3;n-4)=171 - 19 (-4 + n)$, and consequently $\bar{g}(n-4,-3) \leq 0$ if $n \geq 13$.
Furthermore, when $\la \rightarrow -\infty$, $\bar{g}(\la; n-4)\rightarrow +\infty$, which implies $\la'_1 \leq -3$.
Obverse that when $\la \leq -3$,
   $$\bar{g}(\la;n-4)-f(\la;n-5,3)=(-6+n)\la(1+\la)(-2+5\la+2\la^2) > 0.$$
In particular, $f(\la'_1;n-5,3)<0$, which implies $\lamin(\U(n-5,3)^c) < \la'_1$. The result follows. \hfill$\blacksquare$

\begin{lemma}
Given a positive integer $n \ge 20$, for any positive integers $p,q$ such that $p \ge 1,  q \ge 3$ and $p+q=n-2$,
   $$ \lamin(\U(p,q)^c) \ge \lamin (\U( \lceil  (n-2)/2  \rceil, \lfloor (n-2)/2  \rfloor)^c),$$
   with equality \iff  $p= \lceil  (n-2)/2  \rceil$ and $q=\lfloor (n-2)/2  \rfloor$.
\end{lemma}

{\bf Proof:}
Write $\lamin(\U(p,q)^c)=:\la_1$ simply.
By (2.6), we have
  \begin{align*}
  f(\la;p,q)-f(\la;p+1,q-1)&=-\la(2+\la)(-1+2\la)[(p-q+1)(2+\la)+2],\\
  f(\la;p,q)-f(\la;p-1,q+1)&=\la(2+\la)(-1+2\la)[(p-q-1)(2+\la)+2].
  \end{align*}
In addition, $f(-2;p,q)=-10<0$, which implies $\la_1 < -2$.

If $q \geq p+1$, for $\la<-2$ we have
   $$f(\la; p, q ) - f(\la;p + 1, q - 1) > 0.$$
In particular, $f(\la_1;p + 1, q - 1)<0$, which implies $$\lamin(\U(p +1,q - 1)^{c}) < \la_1= \lamin(\U(p,q)^{c}).$$

If $p \geq q+3 (\geq 6)$, by(2.6), we have
  $$f(-3;p,q)= 241- 19p + 23q - 21 pq = 241-19(p-q)+(4-21p)q <0,$$ which implies $\la_1 < -3$.
Observe that
  $f(\la; p, q ) - f(\la;p - 1, q + 1)>0$ when $\la<-3$.
In particular, $f(\la_1;p-1, q+1)<0$, which implies $$\lamin(\U(p-1,q +1)^{c}) < \lamin(\U(p,q)^{c}).$$

To complete the proof, we need to prove $\lamin(\U(p-1,q +1)^{c}) < \lamin(\U(p,q)^{c})$ when $p=q+2$.
In this case, $p=\frac{n}{2}$, $q=\frac{n}{2}-2$, and $$f(\la;p,q)-f(\la;p-1,q+1)=\la(2+\la)(-1+2\la)(4+\la).$$
So it is enough to prove $\la_1 <-4$ or $$f\left(-4;\frac{n}{2}, \frac{n}{2}-2\right)=2376+582n-36n^2<0. $$
If $n \geq 20,$ the above inequality holds, and hence the result follows.
\hfill$\blacksquare$

\section{\bf Main results }
By rearranging the edges of graphs, we first give an maximization of the quadratic form $X^TA(G)X$ among all trees or all unicyclic graphs $G$ of order $n$,
where $X$ is a non-negative or non-positive real vector defined on $G$.

\begin{lemma}
Let $T$ be a tree of order $n$, and let $X$ be a non-negative or non-positive real vector defined on $T$ whose
entries are arranged as $|X_1|\ge |X_2| \ge \cdots \ge |X_n|$ with respect to their moduli.
Then
$$ \sum_{uv \in E(T)} X_uX_v \le \sum_{i=2}^{n} X_1X_i=\sum_{uv \in E(K_{1,n-1})} X_uX_v ,$$
  where $X$ is defined on $K_{1,n-1}$ such that the center has value $X_1$.
If, in addition, $X$ is positive or negative, and $|X_1|>|X_2|$, then the above equality holds only if $T=K_{1,n-1}$.
\end{lemma}

{\bf Proof:}
We may assume $X$ is non-negative; otherwise we consider $-X$.
Let $w$ be a vertex with value $X_1$ given by $X$.
If there exists a vertex $v$  not adjacent of $w$,
   letting $v'$ be the neighbor of $v$ on a path of $T$ connecting $v$ and $w$,
   and deleting the edge $vv'$ and adding a new edge $vu$, we will arrive at a new graph (tree) $T'$, which holds
   $$ \sum_{uv \in E(T)} X_uX_v \le  \sum_{uv \in E(T')} X_uX_v. \eqno(3.1)$$
Repeating the process on the tree $T'$ for the non-neighbors of $w$, and so on,
   we at last arrive at a star $K_{1,n-1}$ with $w$ as its center, and
   $$ \sum_{uv \in E(T)} X_uX_v \le \sum_{uv \in E(K_{1,n-1})} X_uX_v=\sum_{i=2}^{n} X_1X_i. \eqno(3.2)$$

If $X$ is positive, $X_1>X_2$, and $w$ is not adjacent to all other vertices in $T$, then  the inequality (3.1), and hence (3.2),
  cannot hold as an equality. The result follows. \hfill $\blacksquare$

\begin{lemma}
Let $U$ be a unicyclic  graph of order $n$, and let $X$ be a non-negative or non-positive real vector defined on $U$
  whose entries are arranged as $|X_1|\ge |X_2| \ge \cdots \ge |X_n|$ with respect to their moduli.
Then
  $$ \sum_{uv \in E(U)} X_uX_v \le \sum_{i=2}^{n} X_1X_i+X_2X_3=\sum_{uv \in E(S_n^3)} X_uX_v,$$
  where $X$ is defined on $S_n^3$ such that the vertex with degree $n-1$ has value $X_1$,
  and the other two vertices on the triangle have values $X_2,X_3$ respectively.
If, in addition,  $X$ is positive or negative, and $|X_1| > |X_2|$, then the above equality holds  only if $T=S_n^3$.
\end{lemma}

 {\bf Proof:}
We may assume $X$ is non-negative; otherwise we consider $-X$.
Let $w$ be a vertex with value $X_1$ given by $X$.
By a similar discuss to the proof in Lemma 3.1, we have a graph $U'$ of order $n$,
   for which the vertex $w$ is adjacent to all other vertices, and
   $$ \sum_{uv \in E(U)}X_uX_v \le \sum_{uv \in E(U')} X_uX_v=\sum_{i=2}^{n} X_1X_i+X_{u'}X_{v'}, \eqno(3.3)$$
   where $u'v'$ is an edge of $U'$ not incident to $w$.
Surely,
   $$X_{u'}X_{v'} \le X_2X_3. \eqno(3.4)$$
So,
   $$ \sum_{uv \in E(U)}X_uX_v \le \sum_{i=2}^{n} X_1X_i+X_2X_3=\sum_{uv \in E(S_n^3)} X_uX_v. \eqno(3.5)$$

If $X$ is positive, and $X_1 > X_2$, then the equality (3.5) holds only if (3.3) holds,
  which implies $w$ is adjacent to all other vertices  and consequently $U=S_n^3$.
The result follows. \hfill $\blacksquare$

\vspace{3mm}
Next we give a sign property of the first eigenvector of graphs in $\mathscr {U}^{c}_{n}$ for $n \ge 5$.
Note that the first eigenvector of a connected graph with order at least $2$ contains both positive entries and negative entries.

\begin{lemma}
Let $U$ be a unicyclic graph of order $n \ge 5$, and let $X$ be a first eigenvector of $U^c$.
Then $X$ contains at least two positive entries and two negative entries.
\end{lemma}

{\bf Proof:}
Denote by $C_k$ the unique cycle in $U$ of length $k$, and write $\lamin(U^c) =: \la$ simply.
When $n \ge 5$, we know $U^c$ is connected or $U^c=(S_n^3)^c$,
  the latter of which consists of an isolated vertex and a connected non-complete subgraph of order $n-1\ge 4$.
So $X$ contains at least one positive entry and one negative entry.

Assume to the contrary that $X$ contains exactly one positive entry corresponding the vertex ${\bf w}$.
For an arbitrary vertex $u \in N_{U}({\bf w})$, considering the eigenequation (2.2) on $u$ for the graph $U^c$,
  we have
  $$ 0 \le  \lamin(U^c) X_u =\sum_{v \in N_{U^c}(u)} X_v \le 0. $$
So $$X_u=0,  X_v=0,  \hbox{~for each~} u \in N_U({\bf w}), v \in  N_{U^c}(u). \eqno(3.6)$$
If a vertex $r \in N_U(u)$, then
   $N_{U^c}(r)=[ N_{U}(u)\setminus (N_{U}(r) \cup \{r\})] \cup [N_{U^c}(r)\cap N_{U^c}(u)]$,
   and by (3.6),
   $$ \la X_r=\sum_{v \in N_{U^c}(r)} X_v =\sum_{v \in N_{U}(u)\setminus (N_{U}(r) \cup \{r\})} X_v. \eqno(3.7)$$
Similarly, if the vertex $r \notin N_U(u)\cup \{u\}$, then
   $$ \la X_r=\sum_{v \in N_{U^c}(r)} X_v =\sum_{v \in N_{U}(u)\setminus N_{U}(r)} X_v. \eqno(3.8)$$

We divide the discussion into two cases.

{\bf Case 1:}  $k \geq 4$.
As $U$ is a unicycle graph of order $n \geq 5$, there exists a vertex $u \in N_U({\bf w})$ such that $N_U(u) \setminus \{{\bf w}\} \ne \emptyset$.
Let $w' \in  N_{U}(u) \setminus \{w\}$.
Note that $N_U(w') \cap N_U(u)=\emptyset$ as $k \ge 4$
   so that $N_U(w') \subseteq N_{U^c}(u)$, and consequently $X_v=0$ for $v \in N_U(w')$ by (3.6).
By (3.7) we have
   $$ \la X_{w'}=\sum_{v \in N_{U}(u)\setminus (N_{U}(w') \cup \{w'\})} X_v=\sum_{v \in N_{U}(u)\setminus \{w'\}}X_v.$$
By (3.6) and (3.7), we also have
   $$ \la X_{{\bf w}}=\sum_{v \in N_{U}(u)\setminus (N_{U}({\bf w}) \cup \{{\bf w}\})} X_v=\sum_{v \in N_{U}(u)\setminus \{{\bf w}\}}X_v. \eqno(3.9)$$
So
   $$(1+\la) X_{\bf w} =\sum_{v \in N_{U}(u)} X_v=(1+\la)  X_{w'},$$
   which implies $X_{\bf w} = X_{w'}$ as $\la <-1$, a contradiction.

{\bf Case 2:} $k = 3$.
We assert that $w$ has exactly one neighbor in $U$.
Otherwise, let $u,u'$ be two vertices adjacent to ${\bf w}$ in $U$.
As $k=3$, $N_U(u)\setminus \{{\bf w}\} \subseteq N_{U^c}(u')$.
Applying (3.6) on the vertex $u'$, we have $X_v=0$ for each $v \in N_{U^c}(u')$, and hence $X_v=0$ for each $ v \in N_U(u)\setminus \{{\bf w}\}$.
Also by (3.6) $X_u=0$ and $X_v=0$ for each $v \in  N_{U^c}(u)$.
So $X$ contains exactly one nonzero entry, i.e. $X_{\bf w}$; a contradiction.

Let $u$ be the unique neighbor of ${\bf w}$ in the graph $U$.
If there exists a neighbor of $u$ in $U$ other than ${\bf w}$, say $w'$, not lying on the cycle,
   then $N_U(u) \setminus N_U(w')=N_U(u)$, and by (3.7),
   $$\la X_{w'}=\sum_{v \in N_{U}(u)\setminus \{w'\}}X_v.$$
Combining this equality with (3.9), we have $X_{{\bf w}}=X_{w'}$,  a contradiction.
So the vertex $u$ lies on the cycle $C_3$ and has exactly two neighbors except ${\bf w}$.

Now let $w',w''$ be the neighbors of $u$ on the cycle of $U$.
By (3.6) and (3.7), we have
   $$\la X_{w'}= \sum_{v \in N_{U}(u) \setminus \{w',w''\}}X_v = \la X_{w''}, \eqno(3.10)$$
   which implies $X_{w'}=X_{w''}$.
As the graph $U$ has more than $4$ vertices, one of $w',w''$, say $w'$, has a neighbor $\bar{w}'$
   not on the cycle which is also not adjacent to $u$.
By (3.6), $X_{\bar{w}'}=0$, and by (3.8),
   $$ 0= \la X_{\bar{w}'}= \sum_{v \in N_{U}(u) \setminus \{w'\}}X_v. \eqno(3.11)$$
Combining (3.10) with (3.11), we have
    $\la X_{w'}=-X_{w''}=-X_{w'}$, and hence $X_{w'}=X_{w''}=0$ as $\la<-1$.
So, we have $X_v=0$ for each $v \in N_U(u)\setminus \{{\bf w}\}$.
Also by (3.6) $X_u=0$ and $X_v=0$ for each $v \in  N_{U^c}(u)$.
Therefore $X$ contains exactly one nonzero entry, i.e. $X_{\bf w}$; a contradiction.

By above discussion, we get $X$ contains at least two positive entries.
Similarly if considering $-X$, we also get $X$ contains at least two negative entries.
The result follows.   \hfill $\blacksquare$

\vspace{3mm}
We now arrive at the main result of this paper.

\begin{theorem}
Let $U$ be a unicyclic graph of order $n \ge 20$.
Then
$$\lamin(U^c) \ge \lamin(\U( \lceil  (n - 2)/2  \rceil, \lfloor (n - 2)/2  \rfloor)^c),$$
with equality \iff  $U= \U( \lceil  (n - 2)/2  \rceil, \lfloor (n - 2)/2  \rfloor)^c$.
\end{theorem}

 {\bf Proof:}
Let $X$ be the first eigenvector of $U^{c}$ with unit length.
Denote $V_{+} = \{v \in V(U^{c}): X_v > 0, \}$, $V_{-} = \{v \in V(U^{c}): X_v \le 0\}$,
   both being nonempty and containing at least $2$ elements by Lemma 3.3.
Denote by $U_{+}$ (respectively, $U_{-}$)  the subgraph of $U$ induced by
   $V_{+}$ (respectively, $V_{-}$), $E'$ the set of edges between $V_{+}$ and $V_{-}$ in $U$.
Since $U$ is connected, $E' \neq \emptyset$.
Obviously,
$$\sum_{vv' \in E(U)} X_vX_{v'}  =
   \sum_{vv' \in E(U_+)} X_vX_{v'} + \sum_{vv' \in E(U_-)} X_vX_{v'}+
   \sum_{vv' \in E'} X_vX_{v'}.      \eqno(3.12)   $$

First assume $|V_-| \ge 3$. The cycle of $U$ may contain the edges of $E'$, or is contained in one of $U_+, U_-$.
Without loss of generality, we assume that the cycle of $U$ is not contained in $U_+$; otherwise we consider the vector $-X$ instead.
Let $U^*$ be a graph obtained from $U$ by possibly adding some edges within $V^+$ and $V^-$,
such that the subgraph of $U^*$ induced by $V^+$, denoted by $U_+^*$, is a tree, and the subgraph of
$U^*$ induced by $V^-$, denoted by $U_-^*$, is a unicyclic graph.

In the tree $U_+^*$,  choose a vertex, say ${\bf u}$,
with maximum modulus among all vertices of $U_+^*$.
By Lemma 3.1, we will have  a star, say $K_{1,p}$ centered at ${\bf u}$,  where $p+1=|V^+| \ge 2$, which holds
$$ \sum_{vv' \in E(U_+)} X_vX_{v'} \le \sum_{vv' \in E(U_+^*)} X_vX_{v'}
   \le \sum_{vv' \in E(K_{1,p})} X_vX_{v'}. \eqno(3.13)$$

In the unicyclic graph $U^*_-$, choosing a vertex,  say ${\bf w}$,  with maximum modulus.
By Lemma 3.2, we have  a unicyclic graph  $S_{q+1}^3$, where $q+1=|V_-| \ge 3$ and the vertex ${\bf w}$ joins all other vertices of $S_{q+1}^3$,
which holds
  $$ \sum_{vv' \in E(U_-)} X_vX_{v'} \le  \sum_{vv' \in E(U_-^*)} X_vX_{v'} \le \sum_{vv' \in E(S_{q+1}^3)} X_vX_{v'}.\eqno(3.14)$$

Let ${\bf u}',{\bf w}'$ be the vertices of $U_+,U_-$ with minimum modulus among all vertices of $U_+,U_-$, respectively. Then
$$   \sum_{vv' \in E'} X_vX_{v'} \le X_{{\bf u}'}X_{{\bf w}'}, \eqno(3.15)$$
Now by (3.12-3.15), we have
$$  \sum_{vv' \in E(U)} X_vX_{v'}  \le  \sum_{vv' \in E(K_{1,p})} X_vX_{v'}+ \sum_{vv' \in E(S_{q+1}^3)} X_vX_{v'}+X_{{\bf u}'}X_{{\bf w}'}
\eqno(3.16)$$

Since $p \geq 1$, the vertex ${\bf u}'$ can be chosen within the pendent vertices of $K_{1,p}$ by Lemma 3.1.
 If $q \geq 3$, ${\bf w}'$ can be chosen within the pendent vertices of $S_{q+1}^3$ by Lemma 3.2, then from (3.16) we have
$$\frac{1}{2}X^TA(U)X=\sum_{vv' \in E(U)} X_vX_{v'}  \leq  \sum_{vv' \in E(\U(p,q))} X_vX_{v'}=\frac{1}{2}X^TA(\U(p,q))X.\eqno(3.17)$$
and consequently
$$\begin{array} {ccl}
\lamin(U^c) & = & X^TA(U^c)X =X^T (\mathbf{J}-\mathbf{I})X - X^T A(U)X \\
 & \ge & X^T (\mathbf{J}-\mathbf{I})X - X^T A(\U(p,q))X\\
 &= & X^T  A(\U(p,q)^c)X \\
 & \ge & \lamin(\U(p,q)^c).
 \end{array}\eqno(3.18)
$$
If $q = 2$, that is, $S_{q+1}^3 = C_3$,  by a similar discussion,  we have $\lamin(U^c) \geq \lamin(\U'(n-4)^c)$.
By Lemma 2.1, $\lamin(\U'(n-4)^c) > \lamin(\U(n-5,3)^c)$.

Next we consider the case when $|V_-|=2$. In this case the cycle of $U$ cannot lies in $U_-$.
We form a graph $U^\#$ from $U$ possibly by adding some edges within $V^+$ and $V^-$,
such that the subgraph of $U^\#$ induced by $V^+$ is a unicyclic graph, and the subgraph of
$U^\#$ induced by $V^-$ is exactly $K_2$.
Also similar to the discussion for (3.13-3.18), we have
$\lamin(U^c) \ge \lamin(\U(1,n-3))$.
By Lemma 2.2 and the above discussion,
$$ \lamin(U^c) \ge \lamin(\U(p,q)^c) \ge \lamin(\U( \lceil  (n - 2)/2  \rceil, \lfloor (n - 2)/2  \rfloor)^c).\eqno(3.19)$$

Finally we prove the necessity for the equality in the theorem holding.
Assume  the equalities in (3.19) hold.
Then $p = \lceil  (n - 2)/2  \rceil, q= \lfloor (n - 2)/2  \rfloor$, and consequently only the case of $|V_-| \ge 3$ occurs.
Note that $X$ is also a first eigenvector of $\U(p, q)^c$.
Let $\U(p, q)$ have some vertices labeled as in Fig. 2.1, where $v_2={\bf u}, v_3={\bf u}', v_5={\bf w}, v_4={\bf w}'$.

{\bf Assertion 1:} {\it $X$ contains no zero entries; the vertices $v_2={\bf u}$ and  $v_3={\bf u}'$ are respectively the unique ones in $U_+$ with maximum and minimum modulus,
$v_5={\bf w}$ and $v_4={\bf w}'$ are respectively the unique ones in $U_-$ with maximum  and minimum modulus.}
As $X$ is a first eigenvector of $\U(p, q)^c$, from the assumption at the beginning, $X_{v_i}=:X_i>0$ for $i=1,2,3$ and
$X_{v_i}=:X_i \le 0$ for $i=4,5,6,7$.
By (2.5), $$\la_1(X_4-X_7)=-X_3-X_4<0, ~ \la_1(X_6-X_7)=-2X_6,  ~ \la_1(X_5-X_6)=-X_4-(q-3)X_7,$$ which implies that $X_5<X_6<X_7<X_4 \le 0$.
If $X_4=0$, deleting the edges $v_4v_3,v_4v_5$,  and adding two new edges $v_4v_2,v_4v_7$,
  we get a graph $\U(p+1,q-1)$ with same quadratic as $\U(p,q)$ both associated with $X$,
   which implies $\lamin(\U(p,q)) \ge \lamin(\U(p+1,q-1))$, a contradiction by Lemma 2.2.
So $X_4 <0$.
Also by (2.5), $$\la_1(X_1-X_2)>0, ~  \la_1(X_3-X_1)=X_1-X_3-X_4> X_1-X_3,$$ which implies $X_3 < X_1<X_2$.

{\bf Assertion 2:} {\it $U_+=U^*_+=K_{1,p}$, $U_-=U_-^*=S_{q+1}^3$ and $E_1=\{{\bf u}'{\bf w}'\}$, i.e. $U=\U( \lceil  (n - 2)/2  \rceil, \lfloor (n - 2)/2  \rfloor)$.}
By the Assertion 1, retracing the discussion for (3.14-3.15) and applying Lemma 3.1 and Lemma 3.2, we get $U_+=U^*_+=K_{1,p}$, $U_-=U_-^*=S_{q+1}^3$.
From the discussion for (3.15-3.17), also by Assertion 1, $E_1$ consists of exactly one edge and $E_1=\{{\bf u}'{\bf w}'\}$.
\hfill$\blacksquare$

\vspace{3mm}

It was proved in \cite{fan} that $S_n^3$ is the unique minimizing graph in $\mathscr{U}_n$ when $n\ge 6$.
However, when $n \ge 20$, by Theorem 3.4, the graph $\U( \lceil  (n - 2)/2  \rceil, \lfloor (n - 2)/2  \rfloor)^c$ is
the unique  minimizing graph in $\mathscr{U}^c_n$.
So there exists some difference on the least eigenvalue of unicyclic graphs and its complements.

\end{document}